\documentclass[10pt,a4paper]{article}

\usepackage{amsmath}
\usepackage{amssymb}
\usepackage{amsthm}
\usepackage{color}
\usepackage{enumerate}	
\usepackage{graphicx}

\usepackage[T1]{fontenc}% Brzd\k{e}k

\usepackage{titlesec}
\titleformat*{\section}{\normalsize\bfseries}

\usepackage[top=35truemm,bottom=30truemm,left=25truemm,right=25truemm]{geometry}
\makeatletter
\def\mojiparline#1{
    \newcounter{mpl}
    \setcounter{mpl}{#1}
    \@tempdima=\linewidth
    \advance\@tempdima by-\value{mpl}zw
    \addtocounter{mpl}{-1}
    \divide\@tempdima by \value{mpl}
    \advance\kanjiskip by\@tempdima
    \advance\parindent by\@tempdima
}
\makeatother

\theoremstyle{plain}
  \newtheorem{thm}{Theorem}[section]

  \newtheorem{prop}[thm]{Proposition}

\theoremstyle{definition}
  \newtheorem{defn}[thm]{Definition}
  
  \newtheorem{pblm}[thm]{Problem}
\theoremstyle{remark}
  \newtheorem{remark}[thm]{Remark}

\newcommand{\Amp}{{\operatorname{Amp}}}

\newcommand{\CC}{{\mathbb{C}}}

\newcommand{\ext}{\operatorname{ext}}
\newcommand{\Ext}{\operatorname{Ext}}

\newcommand{\Hom}{\operatorname{Hom}}

\newcommand{\NN}{{\mathbb{N}}}

\newcommand{\Num}{{\operatorname{Num}}}
\newcommand{\Ox}{{\mathcal O}_X}

\newcommand{\Pic}{\operatorname{Pic}}
\newcommand{\PP}{{\bf{P}}}

\newcommand{\QQ}{{\mathbb{Q}}}

\newcommand{\RR}{{\mathbb{R}}}

\newcommand{\Spec}{\operatorname{Spec}}

\newcommand{\mf}{{\mathfrak f}}
\newcommand{\sO}{{\mathcal O}}

\begin{document}
%%%%%%%%%%%%%%%%

%\vspace*{5mm}

% 1ページあたり行数の指定
%\linesparpage{48}

\title{Perspective and open problems \\
on birational properties and singularities of \\
moduli scheme of sheaves on surfaces}
\author{Kimiko Yamada
\thanks{Department of Applied Mathematics, 
Faculty of Scienece,
Okayama University of Science, 
Okayama, Japan.
e-mail address:\ k-yamada@ous.ac.jp}}
\date{}

\maketitle

\noindent
{\bf Abstract.} 
%{\small\bf 
For complex projective smooth surface $X$, let $M$
be the coarse moduli scheme  of rank-two stable sheaves with fixed Chern classes.
Grasping the birational structure of $M$, for example its Kodaira dimension, is a fundamental problem.
However, in the case where $\kappa(X)>0$, the study of this problem has not necessarily been active in recent years.
In this article we survey the study of this problem, especially for the case where $\kappa(X)=1$ and $c_1=0$.
We will also survey some research on the structure of singularities of $M$, and
a minimal model program of $M$.
While explaining motivations, we raise several unsolved problems.\\
%}\\[5mm]
%
%Problem \ref{pblm:1st}, Problem \ref{pblm:sing-canonical-or-not},
%Problem \ref{pblm:deg3-part},\ Problem \ref{pblm:direct-kappa}, and Problem \ref{pblm:MMP}.}\\[5mm]

%[5mm]

\section{Overview of singularities and Kodaira dimension of moduli of sheaves on elliptic surfaces with $c_1=0$}

For a complex projective smooth surface $X$ and an ample line bundle $H$ on it,
there is the coarse moduli scheme $M=M_H({\bf v})$ of rank-two $H$-stable sheaves with fixed Chern classes 
${\bf v}=(2,c_1,c_2)$, and
the coarse moduli scheme $\bar{M}=\bar{M}_H({\bf v})$ of $S$-equivalence classes of 
rank-two $H$-semistable sheaves with fixed Chern classes ${\bf v}$ by \cite{Gi:moduli}. 

\begin{pblm}\label{pblm:1st}
(\cite[Question.1.1]{Yamada:Ellipt})
How is the birational structure of $M$, for example its Kodaira dimension $\kappa(M)$?
\end{pblm}

This problem has been actively studied when $\kappa(X)=-\infty,\ 0$.
When $\kappa(X)=2$, J. Li's work \cite{Li:kodaira} is well known.
\cite[Section 11]{HL:text} is an excellent summary for the birational structure of $M$.\par
In this section, we consider the case where $\kappa(X)=1$ and $X$ is a minimal surface.
In this case, $X$ is an elliptic surface, that is, there is a fibration morphism
$\pi:X\rightarrow C$ to a curve $C$ such that its general fibers are elliptic.
We denote the number of multiple fibers of $\pi$ by $\Lambda(X)$, $d =\chi(\Ox)$, and
the fiber class of $\pi$ by $\mf \in \operatorname{Num}(X)$.
For ${\bf v}=(2,c_1,c_2)$,
we suppose that $H$ is ${\bf v}$-suitable (\cite[Definition 3.1.]{Fr:so3ellipt}).
Roughly speaking, a ${\bf v}$-suitable ample line bundle $H$ is not separated from 
$\mf$ by any wall of type ${\bf v}$ (\cite[Definition 2.1.]{Fr:so3ellipt}) in the nef cone of $X$.
When $c_1\cdot \mf$ is odd, Friedman's work is well known (\cite[Theorem 3.14]{Fr:so3ellipt}, 
\cite{yoshioka1999some}).
It states that
$M_H({\bf v})$ is birationally equivalent to a symmetric product of the Jacobian surface $J^{e+1}(X)$
when $X$ is simply connected elliptic surface with at most two multiple fibers. \par
From here we further assume that $c_1=0$.
The reasons for the major differences between these two cases are as follows.
In this case $C=\PP^1$, and we have the generic fiber $X_{\overline{\eta}}$
over $\overline{\eta}=\operatorname{Spec}\overline{k(\PP^1)}$. Let $E$ be an $H$-stable sheaf on $X$.
It induces a vector bundle $E_{\overline{\eta}}$ on $X_{\overline{\eta}}$.
Thanks to the suitability of $H$, $E_{\overline{\eta}}$ is stable if $c_1\cdot \mf$ is odd. 
On the other hand, $E_{\overline{\eta}}$ is semi-stable but not stable if $c_1=0$.\par
Here, we recall a few words.
Let $V$ be a projective normal variety such that $K_V$ is $\QQ$-Cartier.
The $K$-dimension $\kappa(K,V)$ is defined to be
\[ \max\{\dim\Phi_{|mK_V|}:V \dashrightarrow |mK_V| \bigm| 
 m\in \NN, mK_V \text{ is Cartier and } h^0(mK_V)\neq 0\}.\]
The Kodaira dimension $\kappa(V)$ is
$\kappa(K_{\tilde{V}},\tilde{V})$, where $\tilde{V}$ is a smooth complete variety
birationally equivalent to $V$. $\kappa(V)$ is a birational invariant.
If $V$ has only canonical singularities, then $\kappa(K,V)$ and $\kappa(V)$ are equal.
For the definition of canonical singularities, refer to \cite[Def. 6.2.4.]{Ishii:text} for example.\par
Now we denote $M_H((2,0,c_2))=M(c_2)=M$.
Friedman showed the following.
\begin{thm}(\cite[Sect. 7]{Fri:vb-reg-ellipt}\label{Fri:psi})
Suppose $X$ is generic and $c_2>\max(2(1+p_g), 2p_g(X)+(2/3)\Lambda(X))$.
Then there is such a dense open set $M_0$ of $M(c_2)$ as follows.
$M_0$ is contained in
the good locus $M_{gd}$ of $M(c_2)$ defined by
$$ M_{gd}=\left\{ [E]\in M(c_2) \bigm| \ext^2(E,E)^0= 0 \right\}.$$
There is a morphism $\psi:M_0 \rightarrow \PP^{2c_2-2p_g-1}$ such that
the general fiber of $\psi$ is isomorphic to $\Pic^0(C)$, where $C$ is some curve.
\end{thm}
The map $\psi$ in Theorem \ref{Fri:psi} and the prulicanonical maps of $M(c_2)$
are pretty similar:
\begin{prop}(\cite[Cor.3.5.]{Yamada:Ellipt})\label{prop:psi-Phi}
$\psi$ and the pluricanonical map $\Phi:M(c_2)\dashrightarrow \PP^N$ are
coincident up to a quasi-finite map.
In particular, $\dim(\psi)=2c_2-2p_g-1=\{\dim M+1\}/2$ equals the $K$-dimension $\kappa(K,\bar{M}(c_2))$.
\end{prop}
Thus, we can know the Kodaira dimension $\kappa(M(c_2))$ if $M(c_2)$ is projective
(for example, when $c_2$ is odd) and if $M(c_2)$ admits only canonical singularities.
\begin{pblm}\label{pblm:sing-canonical-or-not}
Are all singularities of $M(c_2)$ canonical or not?
\end{pblm}
Here we recall a classical and fundamental fact
in the deformation theory of sheaves from \cite{Lau:Massey}.
If $E$ is a singular point of $M(c_2)$, then for
$b=\dim\Hom(E,E(K_X))^{\circ}\neq 0$ and $D=\dim\Ext^1(E,E)-\dim\Hom(E,E(K_X))^{\circ}$, 
which is the expected dimension of $M(c_2)$, we have
\begin{equation}\label{eq:formal-moduli-intro}
 \sO_{M,E}^{\wedge}\simeq \CC[[t_1,\cdots, t_{D+b}]]/( F_1,\dots, F_b ),
\end{equation}
where $F_i$ is a power series starting from degree-two terms.\par
Any sheaf $E$ in $M(c_2)$ induces a rank-two vector bundle $E_{\eta}$ with degree $0$
on the generic fiber $X_{\eta}$, where $\eta=\Spec(k(\PP^1))$.
From \cite[Fact 2.11.]{Fri:vb-reg-ellipt},
this can be classified into three cases: \par
\begin{enumerate}[{Case} I]
\item : $E_{\eta}$ has no sub line bundle with fiber degree $0$. 
%\textcolor[gray]{0.5}{(Case (A-1))}
\item :
$E_{\eta}$ has a sub line bundle with fiber degree $0$, but $E_{\eta}$ is not decomposable. 
% \textcolor[gray]{0.5}{(Case (B))}
\item :
$E_{\eta}$ is decomposable into line bundles with fiber degree $0$. 
\end{enumerate}
In Case I, we have the following theorem. In the proof, we use a sufficient condition
for singularities to be canonical \cite[Theorem 4.1.]{Yamada:Ellipt}.

\begin{thm}\label{thm:caseI-canonical}
(\cite[Thm. 1.3.]{Yamada:Ellipt})
Suppose that
$E$ is a singular point of $M(c_2)$ applying to Case I.
If $7(d+2)/4\geq \Lambda(X)$ or $2\geq \Lambda(X)$, then the following holds: \\
(1) Let $G$ be any nonzero $\CC$-linear combination of $F_1,\cdots F_b$ in \eqref{eq:formal-moduli-intro} and
we indicate $G$ as
\begin{equation}\label{eq:defn-h-R}
 G=t_1^2+\cdots+ t_R^2+O(3),
\end{equation}
where $O(3)$ are terms whose degrees are more than $2$, 
and $R$ is an integer depending on $G$. Then $R\geq 2b+1$. \\
(2) $E$ is a canonical singularity of $M(c_2)$. \\
Moreover, there actually exist singularities meeting these conditions on $M(c_2)$ if $c_2\gg 0$.
\end{thm}

As a result, we can know $\kappa(M(c_2))$ in some situation where the structure of $X$ is rather simple:

\begin{thm}\label{thm:Kod-dim}
(\cite[Thm. 1.6.]{Yamada:Ellipt})
We suppose that every fiber of $X\rightarrow \PP^1$ is irreducible.
Also we suppose that $X$ has just two multiple fibers
with multiplicities $m_1 = 2$ and $m_2 = m\geq 3$,\ and $d = \chi({\mathcal O}_X) = 1$.
Then all singularities of $M$ apply to Case I.
As a result, all singularities of $M$ are always canonical singularities from Theorem \ref{thm:caseI-canonical}.
Thus the Kodaira dimension $\kappa(M)$ equals to
$\kappa(K,M)=2c_2-2p_g-1=(\dim M(c_2)+1)/2$ from Proposition \ref{prop:psi-Phi}.
\end{thm}

On the other hand, there are cases where it is not possible to know whether $M(c_2)$ has only canonical 
singularities or not only from the evaluation of the second-order terms in the defining equation:

\begin{thm}\label{thm:R=1}
(\cite[Thm.1.4.]{Yamada:Ellipt})
There is an example of an elliptic surface X satisfying the following.
For every obstructed sheaf $E$ in $M(c_2)$, $(M(c_2),E)$ is always a hypersurface singularity, and so
\begin{equation}\label{eq:formal-moduli-hypersurf-intro}
  \sO_{M,E}^{\wedge}\simeq \CC[[t_1,\cdots, t_{D+1}]]/( F ),
\end{equation}
where $F=t_1^2+\cdots+ t_R^2+O(3)$.
There actually exist locally-free obstructed stable sheaves of Case II satisfying $R=1$ in \eqref{eq:defn-h-R}
for every $c_2\gg 0$.
In this case, we cannot judge if $(M,E)$ is a canonical singularity or not only from the second-order terms 
in the defining equation from \cite[Example 7.4.2., Proposition 5.3.12.]{Ishii:text}.
\end{thm}

\section{Problems to be solved in the future 
--Singularities and Kodaira dimension--}

Ongoingly, we suppose that $\kappa(X)=1$ and $c_1=0$ in this section.
We posed the problem of finding the Kodaira dimension of $M$ in Problem \ref{pblm:1st}.
The Kodaira dimension of $M$ was obtained by Theorem \ref{thm:Kod-dim} 
because we were able to show by Theorem \ref{thm:caseI-canonical}
that all singularities of $M$ are canonical.
In this case, the problem was settled by the evaluation of the second-order terms in the defining equation of 
the moduli scheme.
On the other hand, from Theorem \ref{thm:R=1}, it may not be known from the second-order 
terms alone whether the singularity of $M$ is canonical or not.
Therefore the following problem can be raised.

\begin{pblm}\label{pblm:deg3-part}
Let $E$ be a singularity of the moduli space of stable sheaves on surfaces.
Suppose that we can not judge whether it is a canonical singularity or not only from the second-order terms in the defining equation.
Is there a way to evaluate third order or higher terms in a defining equation? 
By using this method, can we judge it is a canonical singularity or not?
\end{pblm}

We can examine the degree-two part of the defining equation of the moduli space from 
linear map between Ext-modules of sheaves
(\cite{Lau:Massey}. cf. \cite[Fact 2.2.]{Yamada:Ellipt}).
However, examining terms of third order or higher would be more difficult than that.
Some related keywords include 
Massey product (\cite{Lau:Massey}) and DG algebra $R\operatorname{Hom}^{\bullet}(E,E)$ 
(See e.g. \cite{BZ18:formal-K3},\ \cite[p.3]{Kaw23:deform-noncomm}). \par
In building a theory using moduli spaces, one may successfully avoid facing singularities at the front.
Thereby it seems worth considering not only Problem \ref{pblm:deg3-part}, but also the following problem.
%So let me raise the following problem.

\begin{pblm}\label{pblm:direct-kappa}
Can we find a way to evaluate the Kodaira dimension $\kappa(M)$ 
without necessarily facing singularities of $M$?
\end{pblm}

As for this direction, there is a work by J. Li \cite{Li:kodaira}.
In this paper, it is shown that $M$ is of general type
when the base surface $X$ is of general type, $p_g(X)\neq 0$, $c_2\gg 0$ and 
$\dim(M)$ is an even number.
Can we come up with some anwer for Problem \ref{pblm:direct-kappa}
in the above-mentioned case where $\kappa(X)=1$ and $c_1=0$?

\section{Problems to be solved in the future --Minimal model program--}

Next, for a normal projective variety $V$ such that $K_V$ is $\QQ$-Cartier and
the singularities of $V$ are log-terminal,
the minimal model program (MMP) of $V$ is proposed and researched.
(For the definition of log-terminal singularities, refer to \cite[Definition 5.2.7.]{Ishii:text}.
We remark that canonical singularities are log-terminal.)

Let us describe a very rough idea of what the MMP of $V$ is.
If the Kodaira dimension $\kappa(V)$ is positive,
we can contract an extremal ray of $V$ to get a small contraction or flip.
After performing these birational transformations finite times,
we obtain a minimal model $V'$ of $V$, that is, $K_{V'}$ becomes nef, and all singularities of $V'$ are log-terminal.
For a detailed explanation, see for example \cite{KM:birat} and \cite{Matsuki:UTX}, where Flowchart 3-1-15 is very useful.

\begin{defn}\label{defn:flip}
Let $f:W \rightarrow Y$ be a birational proper morphism such that $K_W$ is $\QQ$-Cartier and $-K_W$ is $f$-ample, 
the codimension of the exceptional set $\operatorname{Ex}(f)$ of $f$ is more than $1$, and the relative Picard
number of $f$ is $1$. 
We say a birational proper morphism $f_+ : W_+ \rightarrow Y$ is a {\it flip} of $f$ if 
(1) $K_{W_+}$ is $\QQ$-Cartier, 
(2) $K_{W_+}$ is $f_+$-ample,
(3) the codimension of the exceptional set $\operatorname{Ex}(f_+)$ is more than $1$, and
(4) the relative Picard number of $f_+$ is $1$.
\end{defn}

From here, we suppose that $\kappa(X)>0$ and that $X$ is a minimal surface,
and consider the moduli of stable sheaves.
In this case, $K_X$ is
contained in the closure of the ample cone $\Amp(X)$ in $\Num(X)\otimes \RR$.
We have the following theorem on the birational structure of the moduli scheme $M(H)$
of $H$-stable sheaves of type ${\bf v}$,
and the moduli scheme $\bar{M}(H)$ of S-equilavence classes of $H$-semistable sheaves
of type ${\bf v}$.

\begin{thm}[\cite{Yam:flip},\ \cite{Yam:flip-II}]\label{thm:flip}
Suppose $c_2$ is sufficiently large.\\
(1) Let $H$ and $H'$ be ample line bundles on $X$ devided by just one ${\bf v}$-wall $W$.
Assume that $K_X$ does not lie in $W$, and that $K_X$ and $H'$ lie in the same side
with respect to $W$.
For $a\in (0,1)$, we can define $a$-stability of sheaves on $X$ using $H$ and $H'$, and
there is
a moduli scheme $M(a)$ of $a$-stable sheaves of type ${\bf v}$, and
a moduli scheme $\bar{M}(a)$ of S-equivalence classes of $a$-semistable sheaves of type ${\bf v}$.
The wall-crossing problem of $a$-stability induces the sequences of flips
in the sence of Definition \ref{defn:flip}
\begin{align*}
M(H)=M(a_0) \dashrightarrow M(a_1) \dashrightarrow \cdots & 
\dashrightarrow M(a_{N-1}) \dashrightarrow M(a_N)=M(H'),\quad \text{and} \\
\bar{M}(H)=\bar{M}(a_0) \dashrightarrow \bar{M}(a_1) \dashrightarrow \cdots & \dashrightarrow
\bar{M}(a_{N-1}) \dashrightarrow \bar{M}(a_N)=\bar{M}(H') 
\end{align*}
in a moduli-theoretic way.\\
(2) When $H$ varies in $\Amp(X)$ and gets closer to $K_X$,
after crossing finitely many ${\bf v}$-walls,
it reaches an ample line bundle $\tilde{H}$ such that no ${\bf v}$-wall devides $K_X$ and $\tilde{H}$.
Then the canonical class $K_{\bar{M}(\tilde{H})}$ of $\bar{M}(\tilde{H})$ becomes nef.
\end{thm}

Thus one can regard this natural process described in a moduli-theoretic way 
as an analogy of the minimal model program of $M(H)$.
However it is unknown whether 
$M(\tilde{H})$ admits only log-terminal singularities in general.
Thereby we should notice that we need to verify that all singularities of $M(\tilde{H})$ are log-terminal,
in order to say that this sequence is a genuine MMP of $M(H)$.
Here Problem \ref{pblm:sing-canonical-or-not} appears again.
From Theorem \ref{thm:flip}, we can raise the following problem.

\begin{pblm}\label{pblm:MMP}
We pick an ample line bundle $H$ on $X$ and move it closer to $K_X$,
then we obtain an analogy of MMP of $M(H)$ by Theorem \ref{thm:flip}.\\
(1) Can we investigate how $M(H)$ is improved by this moduli-theoretic sequence of flips?
Especially, the case where $M(H)$ admits only log-terminal singularities
is more desirable, because the above-mentioned birational maps give
the genuine MMP in this case.
Also it is more desirable that the structure of the starting point $M(H)$ is easy to understand.\\
(2) It can be interesting if one starts $H$ near $K_X$ and gradually moves $H$ away from $K_X$,
and one observes flips occuring from wall-crossing as in Theorem \ref{thm:flip}.
\end{pblm}

\begin{remark}
We remark that flips in Theorem \ref{thm:flip} are also Thaddeus flips,
that are birational transformations appearing from the variation of GIT quotients
and linearizations (\cite{Tha:GITflip}).
Also, we remark that Thaddeus flips are not necessarily flips in the sense of Definition \ref{defn:flip} in general.
\end{remark}

\end{document}